\documentclass[12pt, reqno]{amsart}
\usepackage{amsmath}
\usepackage{amssymb}
\usepackage{epsfig}
\usepackage{graphicx}
\usepackage{color}
\definecolor{shadecolor}{gray}{0.875}
\usepackage{amscd}
\usepackage{comment}
\usepackage{enumitem}
\usepackage{mathrsfs}

\usepackage[colorlinks=true]{hyperref}

\numberwithin{equation}{section}

\input xy
\xyoption{all}

\calclayout
\allowdisplaybreaks[3]

\theoremstyle{plain}
\newtheorem{prop}{Proposition}[section]

\newtheorem{theo}[prop]{Theorem}
\newtheorem{coro}[prop]{Corollary}

\newtheorem{lemm}[prop]{Lemma}
\newtheorem{assu}[prop]{Assumption}

\theoremstyle{definition}
\newtheorem{defi}[prop]{Definition}

\newtheorem{conj}[prop]{Conjecture}
\newtheorem{pconj}[prop]{Provisional Conjecture}

\newtheorem{rema}[prop]{Remark}

\newtheorem{exam}[prop]{Example}

\def\Br{\mathrm{Br}}

\def\Eff{\overline{\mathrm{Eff}}}
\def\Pic{\mathrm{Pic}}

\def\Mor{\mathrm{Mor}}

\def\Nef{\mathrm{Nef}}

\def\Supp{\mathrm{Supp}}

\def\Spec{\mathrm{Spec}}

\def\Pic{\mathrm{Pic}}

\def\Sec{\mathrm{Sec}}

\def\CDiv{\mathrm{CDiv}}

\makeatother
\makeatletter

\author{Brian Lehmann}
\address{Department of Mathematics \\
Boston College  \\
Chestnut Hill, MA \, \, 02467}
\email{lehmannb@bc.edu}

\author{Sho Tanimoto}
\address{Graduate School of Mathematics, Nagoya University, Furocho Chikusa-ku, Nagoya, 464-8602, Japan}
\email{sho.tanimoto@math.nagoya-u.ac.jp}

\title[Geometric Manin's conjecture]{Geometric Manin's conjecture\\ in characteristic $p$}

\begin{document}
\date{\today}

\begin{abstract}
Geometric Manin's conjecture for complex Fano varieties describes the structure of their moduli spaces of curves.  We propose a version of this conjecture in characteristic $p$ and describe its connection to Batyrev--Manin--Peyre--Tschinkel's conjecture over global fields.
This survey was written for a volume of the Summer Research Institute in Algebraic Geometry held at Colorado State University in 2025.
\end{abstract}

\maketitle

\setcounter{tocdepth}{1}
\tableofcontents

\section{Introduction}

Geometric Manin's conjecture predicts the structure of the moduli space of curves on a Fano variety.  This conjecture has its roots in several influential discoveries made approximately 40 years ago: the close relationship between curves and birational geometry (\cite{Mori79, Mori82}), asymptotic formulas describing counts of rational points (\cite{FMT89, BM}), and the topological properties of spaces of maps from Riemann surfaces into projective Fano manifolds (\cite{Segal79}).

In this expository paper, our goals are:
\begin{enumerate}
    \item To give a conjectural description of families of curves on Fano varieties in characteristic $p$ while accounting for pathological examples.
    \item To give a careful formulation of Batyrev--Manin--Peyre--Tschinkel's conjecture over a global function field with a clear explanation of Peyre's constant in this setting.
    \item To briefly explain the recent developments of  \cite{DLTT25} to solve new cases of Batyrev--Manin--Peyre--Tschinkel's conjecture over global function fields.
\end{enumerate}
This paper is intended to supplement the material in the recent book \cite{BGandManin}.  Throughout, by a curve on a projective variety $X$ we mean a morphism $s: C \to X$ from a smooth projective geometrically integral curve $C$. 

\subsection{Historical background: curves on Fano varieties}

Over the past 50 years, rational curves have emerged as an essential tool in the study of complex Fano varieties.  Some highlights are:
\begin{itemize}
    \item Mori's solution of Hartshorne's conjecture using Bend-and-Break and the deformation properties of free curves (\cite{Mori79}).  Since Mori's ground-breaking work, rational curves have continued to play a key role in our understanding of positivity of the tangent bundle.
    \item The use of lines and conics in the classification of smooth Fano threefolds by \cite{Isk77, Iskov79, iskov, MM8182, MM83, MM03}.
    \item The proof of the boundedness of smooth Fano varieties by \cite{KMM92}.  
\end{itemize}

Starting from the influential works \cite{KP01} (for homogeneous spaces) and \cite{HRS04} (for hypersurfaces), attention turned toward the explicit study of the moduli space of rational curves on Fano varieties: classifying the irreducible components, studying their dimension and singularities, etc.   
Some examples are \cite{Thomsen98, KP01, CS09, Testa09, Bou12, Bou16, BK13, RY16, LT17, LT18, LTdPI, LTDPII, BLRT20, BJ22, Okamura24, Okamura25, MTiB, BV16, BW23}.
For a long time there was no general expectation about how the irreducible components of the moduli space should behave, e.g.~when they should have the expected dimension or how many components represent a fixed numerical class. The situation in characteristic $0$ was clarified in the paper \cite{LT17} and extended to curves of arbitrary genus in \cite{LRT23}.

\subsection{Historical background: Manin's conjecture on rational points}

Yuri Manin and his collaborators initiated a program (later known as the {\it Manin program}) seeking a geometric explanation of various asymptotic formulas for the counting functions of rational points on smooth Fano varieties. This program led to {\it Batyrev--Manin--Peyre--Tschinkel's conjecture} (or {\it Manin's conjecture} for short) which predicts the precise asymptotic formula for the counting function of rational points of bounded height on smooth Fano varieties. This cojecture has been developed in a series of works \cite{FMT89, BM, Peyre, BT, Peyre03, Peyre17, LST18, LS24}.
This conjecture has been confirmed for various examples including, but not limited to, various homogeneous spaces using harmonic analysis \cite{FMT89, BT98, CLT02, STBT, ST16}, complete intersections of low degree using the circle method \cite{Birch, Skinner, BHB17, FM17}, and various (singular) del Pezzo surfaces using the universal torsor method \cite{delaB02, dlBBD, DBB11, dlBBP, Deren14, DP19, DP20, BD25}.

There have also been many works developing Manin's conjecture over global function fields.  Versions of this conjecture were first formulated by Bourqui and Peyre; in \cite{Bourqui03, Bourqui11, Bour13} Bourqui studied toric varieties and varieties with simple Cox rings and in \cite{Peyre12} Peyre studied flag varieties.  More recent examples include hypersurfaces \cite{Lee, Yamagishi, BV16, BW23, Sawin24, HL25} and del Pezzo surfaces \cite{Bour13, DLTT25, Glas25, Tanimoto25} as well as many others.

When one counts the number of rational points, it is important to consider the {\it exceptional set} and remove the contribution of the exceptional set from the counting function so that the asymptotic formula reflects the global geometry of the underlying pair of a smooth Fano variety $X$ and an ample divisor $L$ on $X$.  
The geometric aspect of Manin's conjecture, particularly concerning the exceptional set, has been developed in \cite{HTT15, LT16, LTT14, LTRes, Sen17b, LST18, Gao23}.  In \cite{LST18} the authors and Sengupta proposed a conjectural description of the exceptional set in Manin's conjecture; using recent advances in higher dimensional algebraic geometry such as the Minimal Model Program (the MMP for short) as developed by \cite{BCHM}  and Birkar's solution to the BAB conjecture in \cite{birkar16, birkar16b} we verified that this proposed set is indeed a thin set as predicted by Peyre in \cite{Peyre03}.
There is an alternative proposal using the notion of freeness of rational points combined with the so called all heights approach; see \cite{Peyre17} and \cite{Peyre21}.

\subsection{Batyrev's heuristic and extensions}

In his influential lecture notes \cite{Bat88}, Batyrev developed a heuristic for the asymptotic behavior of rational points on Fano varieties over global fields.  The heuristic is based on the geometry of moduli spaces of curves.  Suppose that $X$ is a Fano variety defined over a finite field.  Counting the number of degree $d$ maps $s: \mathbb{P}^{1} \to X$ is equivalent to counting the $\mathbb{F}_{q}$-points on the associated irreducible components of the space of rational curves $\Mor(\mathbb{P}^{1},X)$.  We can estimate these numbers once we know the dimension and number of irreducible components in degree $d$.  If we assume that these spaces have the ``expected'' geometry, we obtain  the familiar predicted asymptotic formula $cq^{d}d^{\rho(X)-1}$ for the counting function in Manin's conjecture. 
This is further pioneered by David Bourqui and Emmanuel Peyre. See, e.g., \cite{Bourqui03, Bourqui11, Peyre12, Bour13}.

The geometric meaning of the leading constant was later clarified by Ellenberg and Venkatesh (\cite{EV05}). They observed that such an asymptotic formula naturally follows from homological stability results via the Grothendieck-Lefschetz trace formula.  Furthermore, over $\mathbb{C}$ Cohen--Jones--Segal's conjecture (\cite{CJS00}) predicts that the stable homology of the moduli space of curves does indeed stabilize to the homology of the space of continuous maps, and such stable homology naturally leads to Peyre's constant.

\cite{BT} was the first paper to recognize the importance of the Minimal Model Program in establishing Batyrev's geometric heuristic.  Building on recent advances in the MMP, \cite{LT17} and \cite{LRT23} formulated a set of conjectures in characteristic $0$ which translate Batyrev's heuristics into a precise description of the irreducible components of the spaces of curves on any Fano variety (or the spaces of sections of any Fano fibration).  Also \cite{LRT23} completed the first major step toward solving these conjectures by computing the exceptional set in the geometric setting.

An alternative approach to Batyrev's heuristic is to work in the Grothendieck ring of varieties.  This motivic approach was pioneered particularly by Bourqui in \cite{Bou09,Bourqui10,Bourqui11b} and also by \cite{CLL16}.  The recent development of the motivic Euler product by \cite{Bilu} has opened the way for new advances; see e.g.~\cite{BDH22, Lois25}.

\medskip

\noindent
{\bf Disclaimer:} Geometric Manin's conjecture in characteristic $p$ and Batyrev--Manin--Peyre--Tschinkel's conjecture over global function fields are still under development.  This survey paper contains many speculative observations and conjectures; these may need further corrections in the future.

\medskip

\noindent
{\bf Acknowledgements:}
The authors thank Yusuke Nakamura for answering our questions and Alexei Skorobogatov for his help regarding Proposition~\ref{assu:beta}. 
The authors also thank Tim Santens for his explanation of \cite{Santens}. The authors thank Enhao Feng, Jakob Glas, Natsume Kitagawa, Dan Loughran, Marta Pieropan, and Alexei Skorobogatov for comments on a draft of the paper.
Finally, the authors would like to thank the referee whose comments significantly improved the exposition of the paper.

Brian Lehmann was supported by Simons Foundation grant Award Number 851129.  Sho Tanimoto was partially supported by JST FOREST program Grant number JPMJFR212Z and by JSPS KAKENHI Grant-in-Aid (B) 23K25764.

\section{Preliminaries}

Throughout we will work over a ground field $k$.  A variety is an integral separated scheme of finite type over the ground field.  When we say that the ground field has characteristic $p$, we will always implicitly mean that $p$ is a prime number.  In this case, we denote by $F^{e}$ the $e$-th iterate of the (absolute) Frobenius map on $X$ which is a homeomorphism that raises all functions to the $p^{e}$-th power.   If the ground field $k$ is finite, then for any $k$-scheme $X$ we denote by $\mathrm{Fr}_{k}$ the geometric Frobenius automorphism of $X_{\overline{k}}$, i.e.~the morphism $X_{\overline{k}} \to X_{\overline{k}}$ induced by base change from the automorphism of $\overline{k}$ that is the inverse of $t \mapsto t^{|k|}$. A component of a scheme means an irreducible component unless otherwise stated.

\subsection{Numerical spaces}

Let $X$ be a normal projective variety over a field.  We denote the abelian group of Cartier divisors by $\CDiv(X)$.  A $\mathbb{Q}$-Cartier divisor is an element of $\CDiv(X) \otimes_{\mathbb{Z}} \mathbb{Q}$.  

Two Cartier divisors $D_{1},D_{2}$ are numerically equivalent if for every closed integral curve $C$ on $X$ we have $D_{1} . C = D_{2} . C$; numerical equivalence is written as $D_{1} \equiv D_{2}$. We define $N^{1}(X)$ to be the $\mathbb{R}$-vector space of $\mathbb R$-Cartier divisors up to numerical equivalence.

A $1$-cycle on $X$ is a formal sum of closed integral curves on $X$.  Two $1$-cycles $\alpha_{1},\alpha_{2}$ are numerically equivalent if for every Cartier divisor $D$ we have $D . \alpha_{1} = D . \alpha_{2}$; numerical equivalence is written as $\alpha_{1} \equiv \alpha_{2}$.
We define $N_{1}(X)$ to be the real-vector space of $1$-cycles up to numerical equivalence.  This is a finite-dimensional vector space and it contains a natural lattice $N_{1}(X)_{\mathbb{Z}}$ generated by classes of curves. 

\subsubsection{Curves and divisors}

 The space $N_1(X)$ contains two notable cones:
\begin{itemize}
    \item $\Eff_{1}(X)$ is the closure of the cone generated by all effective $1$-cycles.
    \item $\Nef_{1}(X)$ is the nef cone, i.e.~the cone of all numerical classes $\alpha$ such that $E . \alpha \geq 0$ for every effective Cartier divisor $E$.
\end{itemize}
Both cones are closed, convex, full-dimensional, and pointed.
We also consider analogous cones in $N^1(X)$:
\begin{itemize}
    \item $\Eff^1(X)$ is the closure of the cone generated by effective Cartier divisors.
    \item $\Nef^1(X)$ is the nef cone, i.e.~the cone of all $\mathbb R$-Cartier divisors $D$ such that $D.C\geq 0$ for every curve $C$ on $X$.
\end{itemize}
Again, those cones are closed, convex, full-dimensional, and pointed.
We should also note that $\Eff^1(X)$ is dual to $\Nef_1(X)$ and $\Nef^1(X)$ is dual to $\Eff_1(X)$.

\subsection{Brauer groups}
Since we will be working exclusively with quasi-projective schemes $X$ over fields, the cohomological and Azumaya Brauer groups coincide; we denote this common group by $\Br(X)$.  The algebraic part of the Brauer group $\Br_{1}(X)$ is the kernel of the map $\Br(X) \to \Br(X_{k^{s}})$ induced by base change to the separable closure of the base field. 

Let $X$ be a smooth projective variety over a global field $k$.
Then its adelic space is given by
\[
X(\mathbb A_k) := \prod_{v \in \Omega_k}X(k_v),
\]
as a topological space with the product topology, where $\Omega_k$ is the set of all places of $k$ and $k_v$ is the completion of $k$ with respect to $v \in \Omega_k$.
One can define the {\it Brauer--Manin set}
\[
X(\mathbb A_k)^{\mathrm{Br}(X)}
\]
which contains the set $X(k)$ of rational points.
When $X$ is geometrically rationally connected, one of Colliot-Th\'el\`ene's conjectures predicts that $X(k)$ is non-empty as soon as the Brauer--Manin set is non-empty. Moreover the conjecture predicts that $X(k)$ is dense in $X(\mathbb A_k)^{\mathrm{Br}(X)}$.
We refer to this conjecture as Colliot-Th\'el\`ene's conjecture in this paper.
Readers interested in this conjecture should consult \cite{CTS21}.

\subsection{Thin sets}

Suppose $X$ is a variety over a global field $k$.  One possible notion of a ``small'' subset $Z \subset X(k)$ is a non-Zariski dense subset.  However, this notion is not sufficiently flexible for working with arithmetic questions.  Instead, we will need the notion of a thin set introduced by Serre in the context of the inverse Galois problem. 

\begin{defi}
    Let $X$ be a projective variety over a field $k$.  A thin map is a morphism of projective varieties $f: Y \to X$ such that
    \begin{enumerate}
        \item $f$ is generically finite onto its image, and
        \item $f$ is not birational.
    \end{enumerate}    
    A thin set in $X(k)$ is any subset of a finite union $\cup_{i=1}^{r} f_{i}(Y_{i}(k))$ where $\{f_{i}: Y_{i} \to X\}$ is a finite collection of thin maps.
\end{defi}

The notion of a thin set is not particularly useful over some fields (e.g.~an algebraically closed field), but over global fields it is an important and useful condition.
See \cite{Serre} for more details on the notion of thin sets.

\subsection{Fujita invariants}

The following invariant plays a central role in Manin's conjecture:

\begin{defi}
Let $X$ be a smooth projective variety defined over $k$ and $L$ be a big and nef $\mathbb Q$-divisor on $X$.
The {\it Fujita invariant} or {\it $a$-invariant} of $(X, L)$ is the following invariant:
\[
a(X, L): = \min \{ t \in \mathbb R \, | \, \text{the numerical class $t[L] + [K_X] \in \overline{\mathrm{Eff}}^1(X)$} \}.
\]
When $L$ is nef but not big, we set $a(X, L) = + \infty$.
By \cite[Proposition 2.7]{HTT15}, this is a birational invariant under pullback via a birational morphism between smooth projective varieties. (See also \cite[Proposition 4.1.3]{BGandManin}.)

In characteristic $0$, \cite{BDPP} shows that when $L$ is big and nef $a(X, L)$ is positive if and only if $X$ is geometrically uniruled. In positive characteristic, the same statement holds by \cite[Theorem 1.6]{Das20}.
\end{defi}

\begin{defi}
   Let $X$ be a smooth projective variety defined over $k$ and $L$ be a big and nef $\mathbb Q$-divisor on $X$.
   Let $f : Y \to X$ be a dominant generically finite morphism from a smooth projective variety $Y$.
   We say $f$ is an {\it $a$-cover} if $a(Y, f^*L) = a(X, L)$. 
\end{defi}

Next we define the $b$-invariants:
\begin{defi}
    Let $X$ be a smooth projective variety defined over $k$ and $L$ be a big and nef $\mathbb Q$-divisor on $X$.
    Assume that $X$ is geometrically uniruled.
    The {\it face of $(k, X, L)$} is defined by
    \[
    \mathsf F(X, L) := \Nef_1(X)\cap \{ \alpha \in N_1(X) \, | \, (a(X, L)L + K_X).\alpha = 0\}.
    \]
    The {\it $b$-invariant of $(k, X, L)$} is defined by
    \[
    b(k, X, L) = \dim \langle \mathsf F(X, L) \rangle,
    \]
    where $\langle \mathsf F(X, L) \rangle \subset N_1(X)$ is the subspace generated by $\mathsf F(X, L)$.
    By \cite[Proposition 2.10]{HTT15} as well as \cite[Proposition 4.1.16]{BGandManin}, this is a birational invariant under pullback via a birational morphism between smooth projective varieties.
\end{defi}

The most standard example of Fujita invariants and $b$-invariants is the following:

\begin{exam}    Let $X$ be a smooth weak Fano variety defined over $k$ and let $L = -K_X$. Then we have $a(X, L) = 1$ and $b(k, X, L) = \dim N_1(X) = \rho(X)$.
\end{exam}

\begin{defi}
    Let $X$ be a smooth Fano variety and let $f : Y \to X$ be an $a$-cover. We say $f$ is {\it face-contracting} if the induced map
    \[
    \mathsf F(Y, f^*L) \to \mathsf F(X, L),
    \]
    is not injective.
\end{defi}
\subsection{Globally $F$-regular varieties}

\cite{Raynaud78} provided the first example of the failure of Kodaira vanishing for a smooth projective variety in characteristic $p$.  This failure leads to many other pathologies.  However, if we impose an $F$-splitting assumption then we can sometimes recover certain consequences of Kodaira vanishing.

\begin{defi}
    Let $X$ be a smooth projective variety over an algebraically closed field of characteristic $p$.  We say that $X$ is $F$-split if the natural map $\mathcal{O}_{X} \to F_{*}\mathcal{O}_{X}$ admits a splitting in the category of $\mathcal{O}_{X}$-modules.
\end{defi}

We will mainly use the following stronger property.  

\begin{defi}
    Let $X$ be a normal variety over an algebraically closed field of characteristic $p$.  We say that $X$ is globally $F$-regular if for every effective Cartier divisor $D$ there is a positive integer $e$ such that the natural map
    \begin{equation*}
    \mathcal{O}_{X} \to F^{e}_{*}\mathcal{O}_{X}(D)
    \end{equation*}
    admits a splitting in the category of $\mathcal{O}_{X}$-modules.
\end{defi}

The next theorems show that globally $F$-regular varieties share many important properties with Fano varieties in characteristic $0$.  Indeed, one of the key advantages of globally $F$-regular varieties is the following vanishing result:

\begin{theo}[{\cite[Theorem 6.8]{SS10}}] \label{theo:vanishing}
Let $k$ be a field of characteristic $p$.  Suppose $X$ is a smooth projective geometrically integral $k$-variety that is geometrically globally F-regular.  Then for any big and nef Cartier divisor $L$ on $X$ we have $H^{i}(X,\mathcal{O}_{X}(K_{X}+L)) = 0$ for every $i>0$.
\end{theo}


\begin{coro} \label{coro:oxvanishing}
    Let $k$ be a field of characteristic $p$.  Suppose $X$ is a smooth geometrically integral Fano $k$-variety that is geometrically globally F-regular.  Then $H^{i}(X,\mathcal{O}_{X}) = 0$ for all $i>0$.
\end{coro}

Many techniques in Manin's conjecture (such as the formulation of Peyre's constant in \cite{Peyre, BT, CLT10}) require such a vanishing condition.  Since this condition does not hold for arbitrary Fano varieties in characteristic $p$ -- \cite{Maddock16} gives the counterexample of regular del Pezzo surfaces over imperfect fields -- we will rely on these vanishing results in our cases of interest.

Finally we have the following conjecture:
\begin{conj}
    Let $k$ be a field of characteristic $p$.  Suppose $X$ is a smooth geometrically integral Fano $k$-variety that is geometrically globally $F$-regular.   Let $f : Y \to X$ be a dominant generically finite morphism from a smooth projective variety.
    Then we have $a(Y, -f^*K_{X}) \leq a(X, -K_{X})$.
\end{conj}

To the best of our knowledge, this conjecture was not stated before.
Note that in characteristic $0$, this easily follows from the ramification formula. 
In characteristic $p$, the statement is false without the assumption of globally $F$-regular and Fano; see \cite[Example 1.14]{BLRT23} and Example~\ref{exam:tanaka} for examples.

\section{Rational curves on Fano varieties in characteristic $p$}

As discussed above, the work of Mori and his collaborators revolutionized the study of rational curves on Fano varieties.  The following definition identifies the rational curves with the best possible deformation-theoretic properties; such rational curves play a key role in the theory.

\begin{defi}
Let $X$ be a smooth projective geometrically integral variety over a field $k$.  We say that a rational curve $s: \mathbb{P}^{1} \to X$ is:
\begin{itemize}
    \item free, if $s^{*}T_{X}$ is nef;
    \item  very free, if $s^{*}T_{X}$ is ample.
\end{itemize}
More generally, suppose $C$ is a smooth projective geometrically integral curve over $k$.  For any $r \geq 0$, we say that $s: C \to X$ is $r$-free if every positive rank quotient of $s^{*}T_{X}$ has slope at least $2g(C) + r$.
\end{defi}

A famous question of Koll\'ar asks whether every smooth Fano variety over an algebraically closed field of characteristic $p$ carries a very free rational curve or equivalently, is separably rationally connected (See \cite[Chapter IV, 1.10.1 Comment]{Kollar} and also \cite[Chapter IV, 1.13.5 Exercise]{Kollar} for a counterexample in singular cases.)  We review known results on this question in the next two subsections.

\subsection{Counterexamples}

Although Koll\'ar's question is still open for smooth Fano varieties, there are mildly singular Fano varieties which do not carry any free rational curves at all.  The first examples were given in \cite{Kollar95}; the following example is a particular case of Koll\'ar's construction presented by \cite[Section 5]{Xu12}.

\begin{exam}
\label{exam:tanaka}
We work over $\overline{\mathbb{F}}_2$.  Let $\phi: X \to \mathbb{P}^{2}$ be the blow-up of the seven $\mathbb{F}_{2}$-points of $\mathbb{P}^{2}$.  Then $X$ is a weak del Pezzo surface of degree $2$.

There are exactly seven $(-2)$-curves on $X$ corresponding to the strict transforms of the seven $\mathbb{F}_{2}$-lines on $\mathbb{P}^{2}$.  The contraction of these $(-2)$-curves yields a birational map $\phi: X \to X'$ to a degree $2$ log del Pezzo surface $X'$.  \cite[Theorem 4.1.(6)]{CT18} shows that the anticanonical linear series defines a purely inseparable degree $2$ finite morphism $g: X' \to \mathbb{P}^{2}$.  
According to \cite[Proposition 2.4]{Ekedahl87}, the morphism $g$ corresponds to a $p$-closed foliation $\mathcal{F} \subset T_{X'}$.  Since $\mathcal{F}$ is a saturated subsheaf and $X'$ is a normal surface, we see that $\mathcal{F}$ is reflexive and thus \cite[Equation (3.6)]{Ekedahl87} shows that $\mathcal{F} \cong \mathcal{O}_{X'}(-2K_{X'})$.  Let $\psi: T_{X'} \to \mathcal{Q}$ denote the cokernel of the inclusion $\mathcal{F} \subset T_{X'}$.  \cite[Corollary 3.4]{Ekedahl87} shows that $\mathcal{Q}$ injects into $g^{*}T_{\mathbb{P}^{2}}$ and that the cokernel is locally free.  We conclude that $\mathcal{Q} \cong g^{*}\mathcal{O}_{\mathbb{P}^{2}}(-1)$.

The quotient $\psi$ obstructs the existence of (very) free rational curves in the smooth locus of $X'$.  Indeed, if there were such a curve $s: \mathbb{P}^{1} \to X'$, then a general deformation would be contained in the locus where $\psi$ is a surjective map of locally free sheaves.  Then $s^{*}\psi$ would define a negative quotient of $s^{*}T_{X'}$, a contradiction.
\end{exam}

The previous example is a specific instance of an interesting class of weak del Pezzo surfaces $X$ described by the following lemma.  The surfaces which satisfy the equivalent conditions are classified explicitly by \cite{KN22}.

\begin{lemm}[{\cite[Theorem 1.1]{BLRT23}}]
Let $X$ be a weak del Pezzo surface over an algebraically closed field of characteristic $p$.  Then the following properties are equivalent:
\begin{itemize}
\item Every element in $|-K_{X}|$ is singular.
\item $X$ admits a dominant family of rational curves with larger than expected dimension.
\end{itemize}
\end{lemm}

\begin{proof}
    \cite[Theorem 1.1]{BLRT23} implies that if there is a dominant family of rational curves with larger than expected dimension, then every member of $|-K_S|$ is singular. Conversely, if every member of $|-K_S|$ is singular, then all of them are rational curves of arithmetic genus $1$. Since this linear system has dimension larger than expected dimension, our assertion follows.
\end{proof}

In fact, in all such examples there is a finite purely inseparable morphism $f: Y \to X$ such that $K_{Y/X}$ is not contained in $\overline{\mathrm{Eff}}^1(Y)$.
(See \cite[Theorem 1.9]{BLRT23} for this claim.)
In many cases (such as the example of \cite{Xu12} discussed above)  there is a birational model that carries no free curves at all; such examples will necessarily be poorly behaved for Manin's conjecture.

\subsection{Positive examples}

There are several classes of Fano varieties over algebraically closed fields for which the existence of (very) free rational curves is known.  In some cases, one needs to include an $F$-splitting assumption to obtain the best behavior of curves.

\begin{exam}
Suppose $X$ is a smooth Fano hypersurface in $\mathbb{P}^{n}$.  It is known that $X$ is separably rationally connected when $X$ is general (\cite{Zhu24}), or even when $X$ is a general Fano complete intersection (\cite{CZ14}).

The case of an arbitrary smooth Fano hypersurface is still open.  \cite{ST19} shows separable rational connectedness of all Fano hypersurfaces of index $\geq 2$ with degree less than the characteristic (and a similar statement for Fano complete intersections). However, in general it can be difficult to find a very free rational curve: \cite{Cheng25} shows that the minimal degree of a very free rational curve on a Fano hypersurface cannot be bounded above by a linear function in the dimension or the degree. There is an alternative approach to the question of Koll\'ar using the circle method. Indeed, \cite[Theorem 1.5]{BW23} bounds the dimension of the locus parametrizing non-free rational curves in the moduli space of rational curves on a Fano hypersurface when the dimension is large compared to the degree of a hypersurface. In particular, this shows the existence of free rational curves in these cases.
\end{exam}

\cite{CS22} proves similar results for complete intersections in certain homogeneous varieties.

\begin{exam} \label{exam:surfaces}
Suppose that $X$ is a smooth del Pezzo surface.  Since $X$ is rational, it is separably rationally connected.   
Furthermore the result of \cite[Theorem 1.1]{BLRT23} mentioned above shows that every dominant family of rational curves on $X$ has the expected dimension.

If we assume that $X$ has degree $\geq 2$, or $X$ has degree $1$ and the characteristic satisfies $p \geq 11$, then \cite[Theorem 1.2 and Remark 1.4]{BLRT23} show that the following conditions are equivalent: 
\begin{itemize}
\item $X$ is $F$-split.
\item Every irreducible component of $\Mor(\mathbb{P}^{1},X)$ representing a nef class will generically parametrize free rational curves.
\end{itemize}
\end{exam}

\begin{exam} \label{exam:threefolds}
Suppose that $X$ is a normal projective threefold of Fano type in characteristic $> 5$.  Then \cite[Theorem 1.5]{GNT19} shows that $X$ is rationally chain connected.  

If $X$ is a smooth projective threefold that is globally $F$-regular in characteristic $\geq 11$ and $X$ admits a Mori contraction to a curve or surface, then $X$ is separably rationally connected by \cite[Theorem 0.2]{GLPSTZ15}.
\end{exam}

As in the previous examples, it is reasonable to expect better behavior for Fano varieties under an $F$-splitting assumption:

\begin{conj}
\label{conj:SRC}
Let $X$ be a globally $F$-regular smooth Fano variety over an algebraically closed field of characteristic $p$.  Then $X$ is separably rationally connected.
\end{conj}

To the best of our knowledge, this is the first time such a conjecture has been stated explicitly.  However, it is a folklore conjecture that globally $F$-regular varieties will be rationally chain connected or rationally connected.  This topic has been studied in various papers, e.g., \cite{GLPSTZ15, GNT19}. There are also more recent developments, e.g., \cite{APTWX}.

In fact, we expect more to be true: every rational ray in the interior of $\Nef_{1}(X)$ should be represented by the class of a very free rational curve. As discussed in \cite[Remark 4.14]{LRT23}, this property is an important ingredient in geometric Manin's conjecture.

\section{Sections of Fano fibrations in characteristic $p$}

Manin's conjecture addresses the behavior of rational points on a Fano variety over a number field.  To obtain the best geometric analogue, we should study rational points on a Fano variety over the function field of a curve (particularly over a finite base field).  We will always pass directly to an integral model and focus on the equivalent problem of understanding sections of a Fano fibration.

Geometric Manin's conjecture is the study of the asymptotic behavior of sections of Fano fibrations.
In this section, we give a precise formulation of this conjecture, focusing on the geometric aspects.  This section is adapted from \cite{LRT23}.  Some of the conjectures in this section have been established in characteristic $0$ using birational geometry, but currently the characteristic $p$ versions seem out of reach.

\begin{defi} \label{defi:goodfibration}
Let $k$ be a field of characteristic $p$.  A {\it good Fano fibration} is a morphism $\pi: \mathcal{X} \to B$ with the following properties:
\begin{enumerate}
\item $\mathcal{X}$ is a smooth projective geometrically integral variety.
\item $B$ is a smooth projective geometrically integral curve. 
\item $\pi$ is flat and $\mathcal{O}_{B} \cong \pi_{*}\mathcal{O}_{\mathcal{X}}$.
\item The generic fiber $X_{\eta}$ is a smooth geometrically integral Fano variety that is geometrically globally $F$-regular. 
\item The generic fiber $X_{\eta}$ is a Mori dream space.

\item $\pi$ admits a section.
\end{enumerate}
\end{defi}

We denote the moduli space of sections of $\pi: \mathcal{X} \to B$ by $\Sec(\mathcal{X}/B)$.  For a curve class $\alpha \in N_{1}(\mathcal{X})$, $\Sec(\mathcal{X}/B, \alpha)$ denotes the finite type subscheme parametrizing sections with numerical class $\alpha$. 

\begin{rema}
    
    For a projective smooth family over an uncountable algebraically closed field, being globally $F$-regular is an open condition. Indeed, let $f : \mathcal Y \to T$ be a projective smooth family and let $\mathcal L$ be an $f$-ample divisor on $\mathcal Y$. Then by \cite[Proposition 5.3(1)]{SS10}, each fiber $Y_t$ is globally $F$-regular if and only if the section ring
\[
R_t := \bigoplus_{m} H^0(Y_t, \mathcal O(m\mathcal L)),
\]
is strongly $F$-regular. Then being strongly $F$-regular is an open condition by \cite[Theorem B]{ST25}. Thus our claim.
\end{rema}

It is natural to wonder whether some of the properties in Definition \ref{defi:goodfibration} can be weakened.  However we will focus only on this case where a version of geometric Manin's conjecture seems most likely to hold.

The following key definition identifies the analogue of a free curve in the relative setting.  Just as with freeness, relative freeness is determined by the positivity of the pullback of the tangent bundle.  Note that a section will automatically be contained in the locus where $T_{\mathcal{X}/B}$ is locally free.

\begin{defi}
Let $\pi: \mathcal{X} \to B$ be a good Fano fibration.  A section $s: B \to \mathcal{X}$ is relatively $r$-free if every positive rank quotient of $s^{*}T_{\mathcal{X}/B}$ has slope at least $2g(B) + r$.
\end{defi}

\subsection{Numerical classes}
Our first task is to describe the set of numerical classes of sections of $\mathcal{\pi}$.  This set will be contained in a translate of the subspace 
\begin{equation*}
V = \{ \alpha \in N_{1}(\mathcal{X}) \, | \, \alpha . \mathcal X_b = 0 \}
\end{equation*}
where $\mathcal X_b$ is a general geometric fiber of $\pi$.  Thus we first focus our attention on $V$.

Note that we have an injective linear map $N_{1}(\mathcal X_{\eta}) \to N_{1}(\mathcal{X})$ whose image is contained $V$.  (This map is dual to the surjective restriction map $N^{1}(\mathcal{X}) \to N^{1}(\mathcal X_{\eta})$.)  From now on we will identify $N_{1}(\mathcal X_{\eta})$ with its image in $V$.

\begin{lemm}[{\cite[Lemma 5.2]{LRT23}}] \label{lemm:vnef}
Let $\pi: \mathcal{X} \to B$ be a good Fano fibration.  Then we have
\begin{equation*}
\Nef_{1}(\mathcal X_{\eta}) = V \cap \Nef_{1}(\mathcal{X}).
\end{equation*}
We define $\Nef_{vert,\mathbb{Z}}$ to be $\Nef_{1}(\mathcal X_{\eta}) \cap V_{\mathbb{Z}}$ where $V_{\mathbb Z} = V \cap N_1(\mathcal X)_{\mathbb Z}$.
\end{lemm}
\begin{proof}
    Since we are assuming that $\mathcal X_\eta$ is a Mori dream space, the cone $\overline{\mathrm{Eff}}^1(\mathcal X_\eta)$ is equal to the effective cone of divisors. Then the proof of \cite[Lemma 5.2]{LRT23} applies.
\end{proof}

\begin{rema}
There is a minor subtlety: although $N_{1}(\mathcal{X}_{\eta})_{\mathbb{Z}}$ is contained in $V_{\mathbb{Z}} \cap N_{1}(\mathcal{X}_{\eta})$, the two lattices may not coincide; the difference reflects the monodromy of $\pi$.  In particular, the monoid $\Nef_{vert,\mathbb{Z}}$ may be strictly larger than the monoid $\Nef_{1}(\mathcal X_\eta) \cap N_{1}(\mathcal{X}_{\eta})_{\mathbb{Z}}$.  
\end{rema}

We next turn from $V$ to the possible numerical classes of nef sections (contained in a translate of $V$).  We can put restrictions on the numerical classes of sections as follows.  Suppose $\mathcal X_b$ is a reducible fiber of $\pi$ over a closed point $b$ of $B$.  Every section of $\pi$ must intersect $\mathcal X_b$ at a smooth point, and in particular, will intersect a unique irreducible component of $\mathcal X_b$.

\begin{defi}
Let $\pi: \mathcal{X} \to B$ be a good Fano fibration.  An {\it intersection profile} $\lambda$ consists of a choice of one generically smooth and geometrically integral irreducible component in each fiber of $\pi$.  
The (finite) set of all intersection profiles is denoted by $\Lambda$.
\end{defi}

Alternatively, we can identify an intersection profile with the affine subspace of $N_{1}(\mathcal{X})$ that consists of classes $\alpha$ which have intersection $1$ against the geometrically irreducible component in each geometric fiber identified by $\lambda$ and intersection $0$ against all other components.  Henceforth we will not distinguish between these two different ways of thinking about intersection profiles.

The following definition summarizes the above discussion by identifying the possible classes of nef sections.

\begin{defi}
We define $\mathsf{S}_{\mathcal{X}/B}$ to be the convex hull of the set of nef $\mathbb{Z}$-curve classes which have intersection $1$ against a general geometric fiber $\mathcal X_b$ of $\pi$.  

For each intersection profile $\lambda$, we let $\mathsf{S}_{\lambda}$ denote the convex hull of the set of nef $\mathbb{Z}$-curve classes lying in the affine subspace of $N_{1}(\mathcal{X})$ corresponding to $\lambda$. Note that such an affine subspace is a translate of $N_{1}(\mathcal{X}_{\eta})$.
\end{defi}

\begin{rema}
We emphasize that $\mathsf{S}_{\lambda}$ can be strictly contained in the set of nef $\mathbb{R}$-curve classes which have intersection profile $\lambda$.
\end{rema}

Note that $\mathsf{S}_{\lambda,\mathbb{Z}}:= \mathsf S_\lambda\cap N_1(\mathcal X)_{\mathbb Z}$ is preserved by adding nef curve classes in $V_{\mathbb Z}$.  In particular, Lemma \ref{lemm:vnef} shows that $\mathsf{S}_{\lambda,\mathbb{Z}}$ naturally carries the structure of a module over the monoid $\Nef_{vert,\mathbb{Z}}$.  The following conjecture predicts that this module structure essentially controls $\mathsf S_{\lambda,\mathbb{Z}}$.

\begin{conj}
\label{conj:conetheorem}
Let $\pi: \mathcal{X} \to B$ be a good Fano fibration.  For each intersection profile $\lambda$, $\mathsf{S}_{\lambda}$ is a rational polyhedral convex set whose recession cone is $\Nef_{1}(\mathcal{X}_{\eta})$.
\end{conj}

Over $\mathbb C$, the analogous statement is a direct consequence of \cite[Theorem 5.7]{LRT23}.  (See also \cite[Corollary 5.8]{LRT23} for a slightly different version of Conjecture~\ref{conj:conetheorem} over $\mathbb{C}$.)  \cite[Theorem 5.7 and Corollary 5.8]{LRT23} are proved using the cone theorem of the MMP and the boundedness of singular Fano varieties proved by Birkar (\cite{birkar16b}).  Thus we can expect Conjecture \ref{conj:conetheorem} to hold if these foundational results in the MMP and resolutions of singularities are available.

One consequence of this conjecture is that $\mathsf{S}_{\lambda,\mathbb{Z}}$ is ``sandwiched'' between two translates of $\Nef_{vert,\mathbb{Z}}$.  In other words, suppose we take a sum of exponentials $\sum_{\alpha} Cq^{-K_{\mathcal{X}/B} \cdot \alpha}$ indexed over the lattice points $\alpha \in \mathsf{S}_{\lambda,\mathbb{Z}}$.  If we restrict the index set to a translate of $\Nef_{vert,\mathbb{Z}}$ contained in $\mathsf{S}_{\lambda,\mathbb{Z}}$, the difference in the asymptotic growth rates of the sums will be asymptotically negligible in the sense that it converges to $0$ after being divided by the asymptotic formula predicted by Manin's conjecture.

\subsection{Exceptional set}

It is well-known that in Manin's conjecture one must remove an ``exceptional set'' to obtain the correct count.  
\cite{LT16} and \cite{LST18} predict that the exceptional set comes from certain types of maps $f: \mathcal{Y} \to \mathcal{X}$.  


\begin{defi}
    Let $\pi: \mathcal{X} \to B$ be a good Fano fibration.  Suppose that $\mathcal{Y}$ is a smooth projective $B$-scheme such that $\mathcal{Y}_{\eta}$ is geometrically integral.  A thin $B$-morphism $f: \mathcal{Y} \to \mathcal{X}$ is an {\it exceptional map} if either
  \begin{itemize}
\item $f$ is non-dominant and $a(\mathcal Y_\eta, -f^*K_{\mathcal X/B}|_{\mathcal Y_\eta}) \geq 1$;
\item $f$ is dominant and $a(\mathcal Y_\eta, -f^*K_{\mathcal X/B}|_{\mathcal Y_\eta}) > 1$;
\item $f$ is an $a$-cover and $\kappa(-f^*K_{\mathcal X/B}|_{\mathcal Y_\eta} + K_{\mathcal Y_\eta}) > 0$;
\item $f$ is an $a$-cover with $\kappa(-f^*K_{\mathcal X/B}|_{\mathcal Y_\eta} + K_{\mathcal Y_\eta})=0$ and is geometrically non-Galois, or;
\item $f$ is an $a$-cover with $\kappa(-f^*K_{\mathcal X/B}|_{\mathcal Y_\eta} + K_{\mathcal Y_\eta})=0$, geometrically Galois, and face-contracting.
  \end{itemize}
  Here we say that $f$ is geometrically Galois if $f_{\overline{k}} : \mathcal Y_{\overline{k}} \to \mathcal X_{\overline{k}}$ induces a Galois extension of the function fields.
\end{defi}

\begin{rema}
    We expect that certain pathologies of $\Sec(\mathcal{X}/B)$ -- such as the existence of infinitely many irreducible components of $\Sec(\mathcal{X}/B)$ with too large dimension -- should be accounted for by exceptional maps, in the sense that such irreducible components should be the images of sections on some $\mathcal{Y}$.  See \cite[Theorem 1.3]{LRT23} for a statement of this type over $\mathbb{C}$.
\end{rema}

\begin{defi}
Let $\pi: \mathcal{X} \to B$ be a good Fano fibration.  
For a numerical class $\alpha \in \mathsf{S}_{\mathcal{X}/B}$, an irreducible component $M$ of $\Sec(\mathcal X/B, \alpha)$ is an {\it exceptional component} if there is an exceptional map $f : \mathcal Y \to \mathcal X$
and a component $N$ of $\Sec(\mathcal Y/B)$ such that $f$ induces a dominant map
\[
f_* : N \to M.
\]
A component $M$ is a {\it Manin component} if it is not exceptional.
\end{defi}

\subsection{Manin components} 

The final step is to understand the structure of the set of Manin components $M_{\alpha}$.   
The key task facing us is to identify precisely which classes in $\mathsf{S}_{\mathcal{X}/B}$ represent Manin components and how many Manin components there are in each numerical class.

First of all we expect that Manin components parametrize free sections. 

\begin{conj}\label{conj:exceptionalset}
Let $\pi: \mathcal{X} \to B$ be a good Fano fibration.
For each intersection profile $\lambda$, there is an $\alpha_0 \in \mathsf{S}_{\lambda,\mathbb{Z}}$ such that every Manin component $M_{\alpha}$ representing a class $\alpha \in \alpha_0 + \Nef_{vert,\mathbb{Z}}$ will generically parametrize free sections.
\end{conj}

In characteristic $0$, this conjecture is stated in {\cite[Geometric Manin's conjecture 3]{LRT23}. Next, the following conjecture allows us to count the number of Manin components.

\begin{conj}
\label{conj:freecurves}
Let $\pi: \mathcal{X} \to B$ be a good Fano fibration.  For each intersection profile $\lambda$, there is an $\alpha_0 \in \mathsf{S}_{\lambda,\mathbb{Z}}$ such that every \emph{algebraic} equivalence class of curves contained in $\alpha_0 + \Nef_{vert,\mathbb{Z}}$ is represented by exactly one Manin component. 
\end{conj}

In characteristic $0$, this conjecture is stated in \cite[Geometric Manin's conjecture 4]{LRT23}.

We emphasize that it is algebraic and not numerical equivalence that appears in Conjecture \ref{conj:freecurves}. (See \cite[Theorem 1.1]{Okamura25} for an example demonstrating that the number of Manin components for each numerical class is equal to the number of algebraic classes representing that numerical class.)  Currently the best known results which relate algebraic equivalence and the existence of families of curves are due to \cite{Tian25} and \cite{KT25}; see \cite[Theorem 1.1]{KT25} for a precise result in this direction.  Extending these results to solve Conjecture \ref{conj:freecurves} is an important but challenging problem.


\begin{rema}
The main geometric structure carried by spaces of sections is the breaking-and-gluing structure.  That is, using Bend-and-Break one can break off $\pi$-vertical rational curves from any section that has sufficiently many deformations.  Conversely, starting from any section one can obtain a new section by gluing on sufficiently many $\pi$-vertical free rational curves and smoothing.

In practice, the validity of Conjecture \ref{conj:freecurves} depends on whether the general fibers of $\pi$ admit sufficiently many very free rational curves.  This is one reason we have restricted our attention to fibrations whose generic fiber is globally $F$-regular.
\end{rema}

\section{Manin's conjecture over global function fields}

In this section, we state two versions of Manin's conjecture over global function fields. The first version is the standard one based on ideas of \cite{FMT89, BM, Peyre, BT, Peyre03, Peyre17, LST18, LS24} for number fields. The second version is based on the language of Manin components and the all-height-approach of Peyre \cite{Peyre17}.

\subsection{The standard formulation}

Let us describe the set up:
let $k = \mathbb F_q$ be a finite field and let $B$ be a smooth geometrically integral projective curve defined over $k$ with function field $K(B)$.
Let $\pi : \mathcal X \to B$ be a good Fano fibration.
We are interested in counting rational points of the set
\[
\mathcal X_\eta(K(B)).
\]
To this end, we discuss how to define the correct counting function for this set.  
First, unlike over number fields, the height function only takes the values of $q^n$ where $n$ is an integer. For this reason it is necessary to consider the following definitions: 
\begin{defi}
    We define the {\it minimal degree} by 
    \[
     m(\pi) := \min\{ -K_{\mathcal X/B}.C \, | \, [C] \in \Sec(\mathcal X/B)(k) \}.
    \]
    By the Northcott property as in \cite[Lemma 2.6]{LTdPI} and \cite[Lemma 2.2]{LTDPII}, this minimum  is well-defined.
    Next we define the {\it index} by 
    \[
    r(\pi) := \min \{ -K_{\mathcal X/B}.\alpha \, | \,\alpha \in   N_1(\mathcal X)_{\mathbb Z}, \mathcal X_b.\alpha = 0, -K_{\mathcal X/B}.\alpha > 0\}
    \]
    where $\mathcal X_b$ is the class of a geometric fiber of $\pi$. 
    Note that for any section $[C] \in  \Sec(\mathcal X/B)(k)$, we have
    \[
    -K_{\mathcal X/B}.C = m(\pi) + r(\pi)d,
    \]
    for some non-negative integer $d \geq 0$.
    Indeed, let $C_{\min}$ be a section achieving $m(\pi)$.
    Then $\alpha = [C] - [C_{\min}]$ satisfies 
    \[
    \alpha \in   N_1(\mathcal X)_{\mathbb Z}, \quad \mathcal X_b.\alpha = 0, \quad -K_{\mathcal X/B}.\alpha \geq  0.
    \]
    For $\alpha$ satisfying these conditions the possible intersection numbers against $-K_{\mathcal{X}/B}$ will be the intersection of $\mathbb{N} \cup \{0\}$ with a subgroup $r(\pi)\mathbb Z$ of $\mathbb{Z}$.  Our claim follows.

    Next we fix an intersection profile $\lambda$. 
     We define the {\it minimal degree with respect to $\lambda$} by 
    \begin{align*}
     &m(\pi, \lambda) := \\ &\min\{ -K_{\mathcal X/B}.C \, | \, [C] \in \Sec(\mathcal X/B)(k), \text{ the intersection profile of $C$ is $\lambda$}\}.
    \end{align*}
    Next we define the {\it intersection profile index} 
    by 
    \[
    r(\pi)' := \min \{ -K_{\mathcal X/B}.\alpha \, | \,\alpha \in   N_1(\mathcal X_\eta)\cap N_1(\mathcal X)_{\mathbb Z}, -K_{\mathcal X/B}.\alpha > 0\}.
    \]
    Note that $r(\pi)'$ does not depend on $\lambda$.  For any section $[C] \in  \Sec(\mathcal X/B)(k)$ with intersection profile $\lambda$, we have
    \[
    -K_{\mathcal X/B}.C = m(\pi, \lambda) + r(\pi)'d,
    \]
    for some non-negative integer $d \geq 0$.
    We have $r(\pi) | r(\pi)'$ since
    \[
     N_1(\mathcal X_\eta)\cap N_1(\mathcal X)_{\mathbb Z} \subset N_1(\mathcal X)_{\mathbb Z}\cap\{\mathcal X_b.\alpha = 0\}.
    \]
\end{defi}

We introduce the following definition to specify the exceptional set: 

\begin{defi}[{\cite[Definition 4.7]{LRT23}}]
    Let $f : \mathcal Y_\eta \to \mathcal X_\eta$ be a thin map from a geometrically integral smooth projective  variety over $K(B)$. We say $f$ is a {\it breaking thin map} if either there is a strict inequality
    \[
    (a(\mathcal X_\eta, -K_{\mathcal X_\eta}), b(K(B), \mathcal X_\eta, -K_{\mathcal X_\eta})) < (a(\mathcal Y_\eta, -f^*K_{\mathcal X_\eta}), b(K(B), \mathcal Y_\eta, -f^*K_{\mathcal X_\eta}))
    \]
    in the lexicographic order or if equality holds and $f$ is an exceptional map.   
   
\end{defi}

Inspired by \cite[Section 5]{LST18}, we define the {\it exceptional set} to be
\[
\mathsf Z := \bigcup_{f} f(\mathcal Y_\eta (K(B))) \subset \mathcal X_\eta(K(B)) = \Sec(\mathcal X/B)(k)
\]
where $f$ runs over all breaking thin maps. The following conjecture has been proved in \cite[Theorem 5.7]{LST18} in characteristic $0$: 

\begin{conj}
    The exceptional set $\mathsf Z \subset \mathcal X_\eta(K(B))$ is a thin subset of $\mathcal X_\eta(K(B))$.
\end{conj}

This conjecture was not stated in characteristic $p$ before, but it is natural to expect from the point of view of \cite{Peyre03, LST18}.  

With these definitions, we can set up the counting function whose asymptotic behavior we are hoping to understand:

\begin{defi}
    The standard counting function assigns to any 
    non-negative integer $d$ the quantity
    \[
     N_{\mathrm{stan}}(\pi, d):= \# \{ [C] \in (\Sec(\mathcal X/B)(k) \setminus \mathsf Z) \, | \, -K_{\mathcal X/B}.C \leq m(\pi) + r(\pi)d \}.
    \]
\end{defi}

Now our goal is to describe the prediction of the asymptotic formula for this counting function as $d \to \infty$.  
First, we introduce the alpha constant and beta constant:
\begin{defi}
    Let $\pi : \mathcal X \to B$ be a good Fano fibration.
    To the vector space of real $1$-cycles $N_1(\mathcal X_\eta)$ we assign the Lebesgue measure such that the fundamental domain of the lattice $N_1(\mathcal X_\eta) \cap V_{\mathbb Z}$
    has volume equal to $1$. Let $\mathsf C \subset \mathrm{Nef}_1(\mathcal X_\eta) \subset N_1(\mathcal X_\eta)$ be a closed cone. We define the {\it alpha constant} of $\mathsf C$ as
    \[
    \alpha(\mathcal X_\eta, \mathsf C) = \left( \dim \, N_1(\mathcal X_\eta) \right) . \mathrm{vol}(\{\alpha \in \mathsf C \, | \, -K_{\mathcal X/B}.\alpha \leq 1\}).
    \]
    When $\mathsf C = \mathrm{Nef}_1(\mathcal X_\eta)$, we write $\alpha(\mathcal X_\eta, \mathsf C) = \alpha(\mathcal X_\eta)$.
    \end{defi}

    \begin{defi}
    Let $\pi: \mathcal{X} \to B$ be a good Fano fibration.  The {\it beta constant} is
    \[
    \beta(\mathcal X_\eta) = \#(\mathrm{Br}(\mathcal X_\eta)/\mathrm{Br}(K(B))).
    \]
    
\end{defi}
A natural question is whether the beta constant is finite under the assumption that $\mathcal X_\eta$ is smooth Fano and geometrically globally $F$-regular. The following proposition is an affirmative answer to this question:

\begin{prop}
\label{assu:beta}
Let $\pi: \mathcal{X} \to B$ be a good Fano fibration.  Then $\beta(\mathcal{X}_{\eta})$ is finite.
\end{prop}

\begin{proof} 
    We address the algebraic and transcendental contributions to $\beta(\mathcal{X}_{\eta})$ separately.  
    We first claim that $\Pic(\mathcal{X}_{k(\eta)^{s}})$ is torsion-free.  By \cite[Corollary 5.1.3]{CTS21}, we have $\Pic(\mathcal{X}_{k(\eta)^{s}})\cong \Pic(\mathcal{X}_{\overline{\eta}})$ and so it suffices to show that $\Pic(\mathcal{X}_{\overline{\eta}})$ is torsion free.
   Suppose that $D$ is a divisor on $\mathcal{X}_{\overline{\eta}}$ that is numerically trivial.  By Grothendieck-Riemann-Roch and Corollary \ref{coro:oxvanishing} $\chi(\mathcal{O}_{\mathcal{X}_{\overline{\eta}}}(D)) = \chi(\mathcal{O}_{\mathcal{X}_{\overline{\eta}}}) = 1$.  By Theorem \ref{theo:vanishing} we see that $H^{0}(\mathcal{X}_{\overline{\eta}},\mathcal{O}_{\mathcal{X}_{\overline{\eta}}}(D)) = 1$.  Thus $D \sim 0$.
      By torsion-freeness, \cite[Chapter VII, Proposition 4]{Serrelocal} and \cite[Chapter VIII, Corollary 2]{Serrelocal} show that $H^{1}(K(B),\Pic(\mathcal{X}_{k(\eta)^{s}}))$ is finite.  Since $\mathcal{X}_{\eta}$ has a rational point, the Hochschild-Serre spectral sequence 
    \begin{equation*}
        0 \to \Br(K(B)) \to \Br_{1}(\mathcal{X}_{\eta}) \to H^{1}(K(B),\Pic(\mathcal{X}_{k(\eta)^{s}})) \to 0
    \end{equation*}
 implies the finiteness of $\Br_{1}(\mathcal{X}_{\eta})/\Br(K(B))$.

    To prove that $\beta(\mathcal{X}_{\eta})$ is finite, by \cite[Theorem 5.2.5(ii)]{CTS21} it suffices to prove the finiteness of the geometric Brauer group $\mathrm{Br}(\mathcal X_{\overline{\eta}})$.  
    Since a smooth Fano variety over an algebraically closed field is rationally chain connected, an argument of Starr as explained in \cite[Proposition 4.2]{GJ18} shows that there is some integer $N$ that uniformly bounds the order of any element of $\mathrm{Br}(\mathcal X_{\overline{\eta}})$. 
    \cite[Proposition 5.2.3 and Corollary 5.2.8]{CTS21} show that the the prime to $p$ part of $\mathrm{Br}(\mathcal X_{\overline{\eta}})$ is finite. For the $p$-primary part, it follows from \cite[Appendix Theorem A.1]{Skorobogatov} that we have an isomorphism
    \[
    \mathrm{Br}(\mathcal X_{\overline{\eta}})\{p\} \cong (\mathbb Q_p/\mathbb Z_p)^n\oplus H^3(\mathcal X_{\overline{\eta}}, \mathbb Z_p(1))\{p\},
    \]
    where $M\{p\}$ denotes the $p$-primary part of a group $M$.
    Since $ \mathrm{Br}(\mathcal X_{\overline{\eta}})\{p\}$ is annihilated by some power of $p$, $n$ has to be zero.
    We claim that $$H^3(\mathcal X_{\overline{\eta}}, \mathbb Z_p(1))\{p\}$$
    can be interpreted as the extension of a finite group by the group of $\overline{K(B)}$-points on a connected quasi-algebraic unipotent group $\mathsf U$. Indeed, since we assume that $\mathcal X_{\overline{\eta}}$ is globally $F$-regular, it is also $F$-split. 
    Since the Picard group of $\mathcal X_{\overline{\eta}}$ is torsion free, it follows from \cite[Theorem 7.1]{Joshi} that the second crystalline cohomology is exotic torsion free. \cite[Section II.6.7]{Illusie} identifies the divisorial torsion of the second crystalline cohomology as $P^{1}H^{2}(X/W)_{\mathrm{tors}} \oplus H^{2}(X/W)_{v}$. The first summand vanishes by \cite[Proposition II.6.8]{Illusie} since the N\'eron--Severi group is torsion-free.  The second summand vanishes by \cite[Proposition II.6.6]{Illusie} since Corollary~\ref{coro:oxvanishing} implies the vanishing of $H^{2}(W\mathcal O_{\mathcal X_{\overline{\eta}}})$. 
    Moreover by \cite[Theorem 1.1]{Petrov}, the Hodge-to-de Rham spectral sequence of $\mathcal X_{\overline{\eta}}$ degenerates at the first page. Thus the proof of \cite[Proposition 1.5 (ii)]{GSY} justifies the claim.

    The Picard scheme is smooth by \cite[Corollary 7.2]{Joshi} and the vanishing of Witt cohomology.  Thus \cite[Proposition 3.4]{GSY} shows that the formal Brauer group of $\mathcal{X}_{\overline{\eta}}$ is representable and the dimension of the tangent space is $h^2(\mathcal X_{\overline{\eta}}, \mathcal O_{\mathcal X_{\overline{\eta}}})$.  Then \cite[Equation (19)]{GSY} shows that
    \[
    h^2(\mathcal X_{\overline{\eta}}, \mathcal O_{\mathcal X_{\overline{\eta}}}) \geq \dim \mathsf U.
    \]
    In the notation of \cite[Equation (19)]{GSY} our $\mathsf U$ is $\bold U_{\mathcal X_{\overline{\eta}}}$ and we are deducing the above inequality from
    \[
    \dim \mathsf U = \dim \bold U_{\mathcal X_{\overline{\eta}}} \leq \dim \widehat{\mathrm{Br}}(\mathcal X_{\overline{\eta}})_{\mathrm{red}}\leq h^2(\mathcal X_{\overline{\eta}}, \mathcal O_{\mathcal X_{\overline{\eta}}}).
    \]
    Since we have $h^2(\mathcal X_{\overline{\eta}}, \mathcal O_{\mathcal X_{\overline{\eta}}}) = 0$, $\dim \mathsf U$ is $0$. Thus our assertion follows.
\end{proof}

The main ingredient of the leading constant of the asymptotic formula is the Tamagawa number which was first introduced by Peyre \cite{Peyre} for the anticanonical divisor and by Batyrev--Tschinkel \cite{BT} for arbitrary big divisors. A modern account is \cite{CLT10}, and we closely follow its exposition. First we need to introduce the local Tamagawa measures:
\begin{defi}
    Let $\pi : \mathcal X \to B$ be a good Fano fibration. Let $b\in |B|$ be a closed point and denote the completion of $K(B)$ with respect to $b$ by $K(B)_b$.
    Let $\omega$ be a top degree rational differential form on $\mathcal X_{K(B)_b}$ and $\mathrm{div}(\omega)$ be its corresponding divisor.
    Let $\Omega$ be the flat closure of $\mathrm{div}(\omega)$ in $\pi_b: \mathcal X_{\mathfrak o_b} \to \Spec \, \mathfrak o_b$ where $\mathfrak o_b$ is the ring of integers for the local field $K(B)_b$ and $\pi_b$ is the base change of $\pi$. We define the local function
    \[
    \|\omega\| : (\mathcal X\setminus \Supp(\Omega))(K(B)_b)  \to \mathbb R_{>0},
    \]
    that assigns to any $x \in (\mathcal X\setminus \Supp(\Omega))(K(B)_b)$ the quantity
    \[
    \|\omega\|(x) = q_b^{-v_b(x^*\Omega)},
    \]
    where $x : \Spec \, \mathfrak o_b \to \mathcal X_{\mathfrak o_b}$ is the corresponding integral point, $v_b$ is the discrete valuation of $\mathfrak o_b$ and $q_b$ is the size of the residue field.

    Now $\omega$ induces a Radon measure $|\omega|$ on $\mathcal X(K(B)_b)$ and we define the local Tamagawa measure by
    \[
    \tau_{\mathcal X, b} = \frac{|\omega|}{\|\omega\|}.
    \]
    Note that this definition does not depend on the choice of $\omega$ and it only depends on our model $\pi_b : \mathcal X_{\mathfrak o_b} \to \Spec \, \mathfrak o_b$.

\end{defi}

Next, we introduce the local convergence factors:

\begin{defi}
     Let $\pi : \mathcal X \to B$ be a good Fano fibration. Let $\Sigma \subset |B|$ be a finite set such that $\pi$ is smooth outside of $\Sigma$. We define the {\it local convergence factor} $\lambda_b$ for $b \in |B|$ by 
     \begin{equation*}
\lambda_{b} = \left\{ \begin{array}{cl} \det(1 - q_{b}^{-1}\mathrm{Fr}_{\kappa(b)} )^{-1} & \textrm{if } b \in B \backslash \Sigma \\ 1 & \textrm{if }b \in \Sigma \end{array},  \right.
\end{equation*}
where $\kappa(b)$ is the residue field at $b$ and $\mathrm{Fr}_{\kappa(b)}$ is the geometric Frobenius acting on $\Pic(\mathcal{X}_{\overline{\kappa(b)}}) \otimes \mathbb{Q}$. We also define 
\begin{align*}
&L_{\Sigma, *}(1, \Pic(\mathcal X_{\overline{\eta}})\otimes \mathbb Q) :=  \\
&\lim_{t \to 1-} (1-t)^{\rho(\mathcal X_\eta)} \prod_{b \in |B \setminus \Sigma|}\det(1-q_b^{-1}t^{|b|} \cdot \mathrm{Fr}_{\kappa(b)}  )^{-1},
\end{align*}
where $|b| = [\kappa(b):k]$ and the limit is taken so that $t$ is a real number less than $1$ approaching to $1$.
This limit exists as a positive real number.
(This fact is attributed to A. Weil and one can find a statement in \cite[Theorem 9.16B]{Rosen}. See also \cite[Corollary 2.4]{CLT10}.)
We define the {\it Tamagawa measure} $\tau_{\mathcal X}$ on the adelic space
\[
\mathcal X_\eta(\mathbb A_{K(B)}):= \prod_{b \in |B|} \mathcal X(K(B)_b),
\]
by 
\[
\tau_{\mathcal X}:= L_{\Sigma, *}(1, \Pic(\mathcal X_{\overline{\eta}})\otimes \mathbb Q) \prod_{b \in |B|}\lambda_b^{-1}\tau_{\mathcal X, b}.
\]
Since we have $h^i(\mathcal X_\eta, \mathcal O_{\mathcal{X}_{\eta}}) = 0$ for $i = 1, 2$, using the Weil estimates 
\cite[Section 2]{CLT10} proves that this measure is well-defined as a Borel measure and that $\tau_{\mathcal X}(\mathcal X(\mathbb A_{K(B)}))$ is a finite positive real number as soon as $\mathcal X(\mathbb A_{K(B)})$ is non-empty. We should note that \cite{CLT10} assumes that the ground field has characteristic $0$, but all definitions and proofs are valid in characteristic $p$ as well. A key to our claim is \cite[Theorem 2.5]{CLT10}, and the crucial ingredient of this theorem is Deligne’s proof \cite{Deligne} of the Weil conjecture.

Finally we define the {\it Tamagawa number} $\tau_{\mathcal X}(-\mathcal K_{\mathcal X_\eta})$ by
\[
\tau_{\mathcal X}(-\mathcal K_{\mathcal X_\eta}) := q^{m(\pi) + n(1-g(B))}\int_{\mathcal X_\eta(\mathbb A_{K(B)})^{\mathrm{Br}(\mathcal X_\eta)}} \, \mathrm d \tau_{\mathcal X},
\]
where $n = \dim \, \mathcal X_\eta$.
Note that by Proposition~\ref{assu:beta}, the set $$\mathrm{Br}(\mathcal X_\eta)/\mathrm{Br}(K(B)),$$ is a finite set, so the Brauer--Manin set $$\mathcal X_\eta(\mathbb A_{K(B)})^{\mathrm{Br}(\mathcal X_\eta)} \subset \mathcal X_\eta(\mathbb A_{K(B)}),$$ 
is an open and closed set on the adelic space.
Colliot-Th\'el\`ene's conjecture predicts that 
\[
\mathcal X_\eta(K(B)) \subset \mathcal X_\eta(\mathbb A_{K(B)})^{\mathrm{Br}(\mathcal X_\eta)},
\]
is dense so that the above integral is expected to be an integration over $\overline{\mathcal X_\eta(K(B))}$.
\end{defi}

Finally we state a provisional version of Manin's conjecture over global function fields. This is analogous to a version over number fields developed in the series of works \cite{FMT89, BM, Peyre, BT, Peyre03, Peyre17, LST18, LS24}. To the best of our knowledge, the following precise version was not stated before, but some version of this conjecture was proposed by David Bourqui and Emmanuel Peyre. See \cite{Bourqui03, Bourqui11, Peyre12} for examples.

\begin{pconj}[Standard Manin's conjecture over global function fields]
\label{conj:standardManin}
     Let $\pi : \mathcal X \to B$ be a good Fano fibration.
     Suppose that $\mathcal X_\eta(K(B))$ is not thin.
     Then
    \[
    N_{\mathrm{stan}}(\pi, d) \sim (1-q^{-r(\pi)})^{-1}\alpha(\mathcal X_\eta)\beta(\mathcal X_\eta)\tau_{\mathcal X}(-\mathcal K_{\mathcal X_\eta})q^{r(\pi)d}(dr(\pi))^{\rho(\mathcal X_\eta)-1},
    \]
    as $d \to \infty$ with $d \in \mathbb Z$.
\end{pconj}

\begin{rema}
    We expect that when 
    \begin{equation}
    \label{equation:simplification}
         r(\pi) = r(\pi)',
    \end{equation}
  is true, Provisional Conjecture~\ref{conj:standardManin} holds. Note that (\ref{equation:simplification}) is satisfied when every fiber of $\pi$ is integral. However, when (\ref{equation:simplification}) fails,  the asymptotic formula may have some periodicity observed by \cite[Example 1.1.4]{LL25} in the context of Malle's conjecture, i.e.~
  the limit exists only if we fix the congruence class of $d$ modulo some $n$. \cite[Theorem 1.1]{Santens} suggests that we might need some averaging over degrees.
  This is the main reason why Conjecture~\ref{conj:standardManin} is provisional.
\end{rema}

\subsection{All height approach}
It is also natural to consider the following modification of our counting problem.  For each algebraic curve class $\alpha \in \mathsf{S}_{\mathcal{X}/B}$ consider the quantity
\[
\frac{\#M_\alpha(k)}{q^{-K_{\mathcal X/B}.\alpha + n(1-g(B))}}
\]
where $M_\alpha$ denotes a Manin component associated to $\alpha$. We would like to have an asymptotic formula describing the limit of this quantity as $-K_{\mathcal X/B}.\alpha$ goes to $\infty$. As before, it is important to exclude some subloci of $M_\alpha$ to obtain the correct leading constant.
To this end, we would like to propose the following definition:
\begin{defi}
    A thin set $Z \subset \mathcal X_\eta(K(B))$ is properly constructible if there are finitely many finite thin $B$-maps $f_i : \mathcal Y_i \to \mathcal X$ from geometrically integral projective $B$-varieties $\mathcal Y_i \to B$ such that
    \[
    Z = \bigcup_i f_i(\mathcal Y_{i, \eta}(K(B))).
    \]
\end{defi}

We conjecture the following, and this is new in the literature:
\begin{conj}
\label{conj:constructiblethin}
    The exceptional set $\mathsf Z \subset \mathcal X_\eta(K(B))$ is a properly constructible thin set.
\end{conj}

We should note that this conjecture is still open even in characteristic $0$. Using this property, one can construct an open subscheme $M_\alpha^\circ \subset M_\alpha$ by removing the exceptional set:

\begin{prop}
    \label{prop:properlyconstructiblethinset}
    Let $f_i : \mathcal Y_i \to \mathcal X$ be finitely many finite thin $B$-maps.
    For any extension $k'/k$,
    let $Z_{k'} \subset \mathcal X_\eta(K(B)\otimes_k k') = \Sec(\mathcal X/B)(k')$ be the properly constructible thin set defined by $f_{i, k'}$'s.
    Then for any irreducible component $M \subset \Sec(\mathcal X/B)$, there exists a unique open subscheme $M^\circ \subset M$ such that for any extension $k'/k$ we have
    \[
    M^\circ(k') = M(k') \setminus Z_{k'}.
    \]
\end{prop}

\begin{proof}
    We may assume that $k$ is algebraically closed.
    Then our claim follows from the fact that for any finite $B$-morphism $f : \mathcal Y \to \mathcal X$, the induced map
    \[
    f_* : \overline{M}_{g(B), 0}(\mathcal Y) \to \overline{M}_{g(B), 0}(\mathcal X)  
    \]
is proper. 
Note that since $f$ is finite, for any stable map that is numerically equivalent to a section its composition with $f$ is still a stable map.
This means $f_{*}$ maps a stable map with reducible domain to a stable map with reducible domain.
Thus our assertion follows.
\end{proof}

We next set up the counting problem. Let $\pi : \mathcal X \to B$ be a good Fano fibration. 
Let $$\mathsf Z \subset \Sec(\mathcal X/B)(k),$$ be the exceptional set. By Conjecture~\ref{conj:constructiblethin}, $\mathsf Z$ is defined by finitely many finite thin $B$-maps $\{f_i\}$.
Let $\lambda$ be an intersection profile.
Recall that $\mathsf S_{\lambda}$ denotes the set of real nef classes of section type which are compatible with the intersection profile $\lambda$.
We fix a a basis of $N_1(\mathcal X)$ consisting of $\mathbb{Z}$-curve classes and work with the max norm on $N_{1}(\mathcal{X})$ with respect to this basis.

\begin{defi}
    We let $\ell_{\lambda} : \mathsf S_{\lambda} \to \mathbb R_{\geq 0}$ denote the rational piecewise linear function that measures the distance to the relative boundary of $\mathsf S_{\lambda}$ (that is, the boundary of $\mathsf S_{\lambda}$ considered as a subset of the translate of $N_{1}(\mathcal{X}_{\eta})$ corresponding to $\lambda$).
\end{defi}

Conjecture \ref{conj:freecurves} predicts that there is a class $\alpha_{0} \in \mathsf S_{\lambda,\mathbb{Z}}$ such that every algebraic equivalence class of curves contained in $\alpha_0 + \Nef_{vert,\mathbb{Z}}$ is represented by exactly one Manin component.  Conjecture \ref{conj:conetheorem} implies that there is a $\mathbb{Z}$-curve class $\alpha_{1}$ in $\alpha_{0} + N_{1}(\mathcal{X}_{\eta})$ such that $\mathsf S_{\lambda} \subset \alpha_{1} + \Nef_{vert,\mathbb{Z}}$.  Since
\begin{equation*}
    \alpha_{0} + \Nef_{vert,\mathbb{Z}} \subset \mathsf S_{\lambda} \subset \alpha_{1} + \Nef_{vert,\mathbb{Z}}
\end{equation*}
we see that if $\alpha \in \mathsf S_{\lambda}$ satisfies $\ell_{\lambda}(\alpha) \geq \Vert \alpha_{0} - \alpha_{1} \Vert$ then $\alpha \in \alpha_0 + \Nef_{vert,\mathbb{Z}}$.  In summary, when $\ell_{\lambda}(\alpha)$ is sufficiently large we will henceforth assume that the prediction of Conjecture~\ref{conj:freecurves} applies to $\alpha$.

Conjecture \ref{conj:constructiblethin} predicts that there are finitely many thin $B$-maps $f_{i}: \mathcal{Y}_{i} \to \mathcal{X}$ that account for the exceptional set.  Let $M^\circ_\alpha \subset M_\alpha$ be the open subscheme obtained by applying Proposition \ref{prop:properlyconstructiblethinset} to this finite set of morphisms.
With these setups, we have the following conjecture. There are very similar conjectures proposed by David Bourqui  (see \cite[Question 1.10]{Bourqui11b} and \cite[Definition 2.2]{Bour13}) and Emmanuel Peyre (see \cite[Question 6.27]{Peyre21}):

\begin{conj}[All height approach version of Manin's conjecture]
\label{conj:allheight}
    There exist constants $c(\mathcal X, \lambda) \geq 0$ and $\delta > 0$ such that we have
    \[
    \#M_\alpha^\circ(k) = c(\mathcal X, \lambda) q^{-K_{\mathcal X/B}.\alpha + n(1-g(B))} + O(q^{-K_{\mathcal X/B}.\alpha -\delta\ell_{\lambda}(\alpha)}),
    \]
    as $\ell_{\lambda}(\alpha) \to \infty$. Here $n = \dim \, \mathcal X_\eta$.
\end{conj}

Moreover, one can provide a precise prediction for the value of $c(\mathcal X, \lambda)$.
Let $\pi : \mathcal X \to B$ be a good Fano fibration. Let $\Sigma_\pi \subset |B|$ be the set of closed points $b \in |B|$ such that the fiber $\mathcal X_b \to \Spec \, \kappa(b)$ is not smooth.
We denote the smooth locus of this fiber by $\mathcal X_b^{\mathrm{sm}}$.
Since for any $B$-morphism $\sigma : \Spec \, K(B)_b \to \mathcal X$, its specialization $\sigma(b)$ lands in $\mathcal X_b^{\mathrm{sm}}(\kappa(b))$, we have the following continuous map:
\[
\Phi : \mathcal X_\eta(\mathbb A_{K(B)})^{\mathrm{Br}(\mathcal X_\eta)} \subset \prod_{b \in |B|}\mathcal X_\eta(K(B)_b) \to \prod_{b \in \Sigma_\pi} \mathcal X_b^{\mathrm{sm}}(\kappa(b)).
\]
For an intersection profile $\lambda$ and $b \in \Sigma_\pi$, let $\mathcal X_{b, \lambda}^{\mathrm{sm}}$ be the component of $\mathcal X_b^{\mathrm{sm}}$ specified by the intersection profile $\lambda$.
Then we define 
\[
\mathcal X_\eta(\mathbb A_{K(B)})^{\mathrm{Br}(\mathcal X_\eta)}_{\lambda} := \Phi^{-1}\left(\prod_{b\in \Sigma_\pi}\mathcal X_{b, \lambda}^{\mathrm{sm}}(\kappa(b))\right),
\]
which is an open and closed subset of $\mathcal X_\eta(\mathbb A_{K(B)})$ because of Proposition~\ref{assu:beta}.
Here is a conjectural description of $c(\mathcal X, \lambda)$:

\begin{conj}
\label{conj:leading}
    
   In Conjecture~\ref{conj:allheight}, we have
    \[
    c(\mathcal X, \lambda) = \frac{\beta(\mathcal X_\eta)}{\mathsf B}\int_{\mathcal X_\eta(\mathbb A_{K(B)})^{\mathrm{Br}(\mathcal X_\eta)}_{\lambda}}\, \mathrm d\tau_{\mathcal X},
    \]
    where $\mathsf B$ is the number of algebraic equivalence classes representing one single numerical class.
    
\end{conj}
This is based on a conjectural description of the leading constant by Peyre (\cite{Peyre}).

\begin{rema}
It is natural to wonder whether $\mathsf{B}$ is equal to $\#\mathrm{Br}(\mathcal X)$. Over $\mathbb C$, \cite[Theorem 10.17]{Voi03} implies that $H^{2, 0}(\mathcal X) = 0$. Thus the exponential sequence shows that 
\[
H^3(\mathcal X, \mathbb Z)_{\mathrm{tors}} \cong \mathrm{Br}(\mathcal X).
\]
In particular, this proves that $\mathrm{Br}(\mathcal X)$ is finite. 
Then by the Poinc\'are duality we have
\[
\#H_2(\mathcal X, \mathbb Z)_{\mathrm{tors}} = \#\mathrm{Br}(\mathcal X).
\]
For a Fano fibration $\mathcal X$, it is expected that the Griffith group of $\mathcal X$ is trivial and the integral Hodge conjecture for curves holds. (This is known in dimension $\leq 3$ by \cite[Theorem 1]{BS83} and \cite[Theorem 2]{Voisin06}.) Thus we are expecting that the constant $\mathsf B$ coincides with the size of $\mathrm{Br}(\mathcal X)$. See \cite[Remark 4.13]{LRT23} for more details.
\end{rema}

\subsection{Relationship between conjectures}
Finally we give a heuristic argument explaining how the all height approach (when combined with earlier conjectures) implies a version of Provisional Conjecture~\ref{conj:standardManin}.

\begin{assu}
    In this subsection we assume the validity of Conjectures~\ref{conj:conetheorem}, \ref{conj:exceptionalset}, \ref{conj:freecurves}, \ref{conj:allheight}, and \ref{conj:leading}.
    We also assume that (\ref{equation:simplification}) holds to avoid some lattice issues.
    \end{assu}

Consider the following counting function:
\[
N_{\mathrm{Manin}}(\pi, d) := \sum_\lambda \sum_{\substack{\alpha \in \mathsf S_{\lambda, \mathbb Z} \\ -K_{\mathcal X/B}.\alpha \leq m(\pi) + dr(\pi)}} \sum_{M_\alpha} \#M^\circ_\alpha(k),
\]
where $\lambda$ runs over all intersection profiles and $M_\alpha$ runs over all Manin components associated to the numerical class $\alpha$.
Note that for distinct intersection profiles $\lambda \neq \lambda'$, we have $\mathsf S_\lambda \cap \mathsf S_{\lambda'} = \emptyset$, so it is justified to count each intersection profile separately.
Moreover $N_{\mathrm{Manin}}(\pi, d) \leq N_{\mathrm{stan}}(\pi, d)$, but they need not agree when there are exceptional components that do not come from breaking thin maps.
We expect that the difference between $N_{\mathrm{stan}}(\pi, d)$ and $N_{\mathrm{Manin}}(\pi, d)$ is asymptotically negligible, in the sense that it converges to $0$ when divided by the asymptotic formula predicted by Manin's conjecture, but do not have a rigorous proof.  
Here we focus on $N_{\mathrm{Manin}}(\pi, d)$.

Let $\epsilon > 0$ be a sufficiently small rational number.
Then we define the following shrunken set:
\[
\mathsf S_{\lambda, \epsilon} := \{ \alpha \in \mathsf S_{\lambda}\, | \, \ell_{\lambda}(\alpha) \geq -\epsilon K_{\mathcal X/B}.\alpha \}.
\]
In the view of Conjecture~\ref{conj:conetheorem}, this is a rational polyhedral convex set with the recession cone
\begin{equation*}
    \Nef_{1}(\mathcal X_\eta)_\epsilon = \{ \alpha \in \Nef_{1}(\mathcal{X}_{\eta}) \, | \, \ell(\alpha) \geq -\epsilon K_{\mathcal X/B}.\alpha \},
\end{equation*}
where $\ell$ is the distance function to the relative boundary of $\Nef_1(\mathcal X_\eta)$ inside $N_1(\mathcal X_\eta)$.
Let $\mathsf S_{\lambda, \epsilon, \mathbb Z} = \mathsf S_{\lambda,\epsilon} \cap N_{1}(\mathcal{X})_{\mathbb{Z}}$.  When $\alpha \in \mathsf S_{\lambda, \epsilon,\mathbb{Z}}$ and $-K_{\mathcal X/B}.\alpha$ is sufficiently large, $\ell_{\lambda}(\alpha)$ will also be sufficiently large so that every algebraic equivalence class contained in $\alpha$ is represented by a unique Manin component.  We denote those Manin components by $M_{\alpha, i}$ for $i = 1, \cdots, \mathsf B$.  Since there are only finitely many classes in $\mathsf S_{\lambda, \epsilon, \mathbb Z}$ to which the above discussion does not apply, when computing asymptotic behavior we may assume that every class $\alpha$ in $\mathsf S_{\lambda, \epsilon, \mathbb Z}$ is represented by $\mathsf B$ Manin components.

Let us consider the counting function for classes in $\mathsf{S}_{\lambda, \epsilon, \mathbb{Z}}$:
\[
\sum_{\substack{\alpha \in \mathsf S_{\lambda, \epsilon, \mathbb Z} \\ -K_{\mathcal X/B}.\alpha \leq m(\pi) + dr(\pi)}} \sum_{i = 1}^{\mathsf B} \#M^\circ_{\alpha, i}(k).
\]
By Conjectures~\ref{conj:allheight} and \ref{conj:leading}, this is equal to 
\begin{align*}
&\sum_{\substack{\alpha \in \mathsf S_{\lambda, \epsilon, \mathbb Z} \\ -K_{\mathcal X/B}.\alpha \leq m(\pi) + dr(\pi)}} \\ &\left(\beta(\mathcal X_\eta)\tau_{\mathcal X}(\mathcal X_\eta(\mathbb A_{K(B)})^{\mathrm{Br}(\mathcal X_\eta)}_{\lambda})q^{-K_{\mathcal X/B}.\alpha + n(1-g(B))} + O(q^{-(1-\delta\epsilon)K_{\mathcal X/B}.\alpha})
\right).
\end{align*}
Then using the counting arguments of \cite[Proposition 4.3]{LRT23} (combined with Conjecture~\ref{conj:conetheorem} and the argument of \cite[Theorem 5.7]{LRT23}), the asymptotic formula for the sum is given by
\[
(1-q^{-r(\pi)})^{-1}\alpha(\mathcal X_\eta, \Nef_1(\mathcal X_\eta)_\epsilon)\beta(\mathcal X_\eta)\tau_{\mathcal X}(\mathcal X_\eta(\mathbb A_{K(B)})^{\mathrm{Br}(\mathcal X_\eta)}_{\lambda})q^{r(\pi)d}(dr(\pi))^{\rho(\mathcal X_\eta)-1},
\]
as $d \to \infty$ with $d \in \mathbb{Z}$. Hence after taking the summation over $\lambda$, we obtain
\[
(1-q^{-r(\pi)})^{-1}\alpha(\mathcal X_\eta, \Nef_1(\mathcal X_\eta)_\epsilon)\beta(\mathcal X_\eta)\tau_{\mathcal X}(-\mathcal K_{\mathcal X_\eta})q^{r(\pi)d}(dr(\pi))^{\rho(\mathcal X_\eta)-1}.
\]
To finish the argument, we introduce one more conjecture which is expected among the experts. However, we could not locate a reference explicitly stating this conjecture:
\begin{conj}
  There exists uniform constants $\mathsf B'$ and $\mathsf C$ such that for any $\alpha \in \mathsf S_{\lambda, \mathbb Z}$,
  the number of Manin components of the class $\alpha$ is bounded by $\mathsf B'$ and for such a Manin component $M_\alpha$, we have
  \[
  \#M_\alpha^\circ(k) \leq \mathsf Cq^{-K_{\mathcal X/B}.\alpha + n(1-g(B))}.
  \]
\end{conj}
Since $\mathsf{S}_{\lambda,\epsilon,\mathbb{Z}} \subset \mathsf{S}_{\lambda, \mathbb{Z}}$ and the difference is controlled by the region $\Nef_{1}(\mathcal{X}_{\eta}) \backslash \Nef_{1}(\mathcal{X}_{\eta})_{\epsilon}$, we have
\begin{align*}
&(1-q^{-r(\pi)})^{-1}\alpha(\mathcal X_\eta, \Nef_1(\mathcal X_\eta)_\epsilon)\beta(\mathcal X_\eta)\tau_{\mathcal X}(-\mathcal K_{\mathcal X_\eta})\\
&\leq \liminf_{d \to \infty}\frac{N_{\mathrm{Manin}}(\pi, d)}{q^{r(\pi)d}(dr(\pi))^{\rho(\mathcal X_\eta)-1}} \\
& \leq \limsup_{d \to \infty}\frac{N_{\mathrm{Manin}}(\pi, d)}{q^{r(\pi)d}(dr(\pi))^{\rho(\mathcal X_\eta)-1}}  \\
&\leq(1-q^{-r(\pi)})^{-1}\alpha(\mathcal X_\eta, \Nef_1(\mathcal X_\eta)_\epsilon)\beta(\mathcal X_\eta)\tau_{\mathcal X}(-\mathcal K_{\mathcal X_\eta})\\
&+(1-q^{-r(\pi)})^{-1}\mathsf C\mathsf B'\alpha(\mathcal X_\eta, \Nef_1(\mathcal X_\eta) \setminus \Nef_1(\mathcal X_\eta)_\epsilon).
\end{align*}
As $\epsilon \to 0$, we conclude
\[
\lim_{d \to \infty}\frac{N_{\mathrm{Manin}}(\pi, d)}{q^{r(\pi)d}(dr(\pi))^{\rho(\mathcal X_\eta)-1}} = (1-q^{-r(\pi)})^{-1}\alpha(\mathcal X_\eta)\beta(\mathcal X_\eta)\tau_{\mathcal X}(-\mathcal K_{\mathcal X_\eta}).
\]

\section{Examples}

Finally, we discuss a few examples where (weaker forms of) the conjectures of the previous sections have been verified.  We do not attempt to give a complete account of the literature, instead we focus on a few notable examples.

\subsection{Toric varieties}

Manin's conjecture for smooth projective toric varieties over $\mathbb{F}_{q}(t)$ has been settled by David Bourqui in \cite{Bourqui03} and \cite{Bourqui11}. The first paper used the method of universal torsors developed by Salberger \cite{Sal98} over number fields 
and the second paper used harmonic analysis on tori following ideas of Batyrev--Tschinkel in \cite{BT96b, BT98}. Indeed, using these methods, Bourqui obtained a meromorphic continuation of the height zeta function.
Then Conjecture~\ref{conj:standardManin} follows from the Tauberian theorem as in \cite[Corollary A.13]{CLT10}.  Bourqui also proves the all height version of Manin's Conjecture in \cite[Section 2.9]{Bourqui11b}.

Bourqui's pioneering work has been quite influential.  In fact \cite{Bou09} proved a motivic version of Manin's Conjecture for toric varieties; his work has been revisited in \cite{BDH22} and \cite{Lois25}.

\subsection{Low degree hypersurfaces}
Manin's conjecture over global function fields has been proved for low degree hypersurfaces in $\mathbb P^n$ with $p$ being greater than the degree using the Hardy--Littlewood circle method which is a technique from analytic number theory. See \cite{Lee, Yamagishi, BV16, BW23, Sawin24, HL25} for more details.

\subsection{Homological sieve method}

Let $k = \mathbb F_q$ be a finite field and $S$ be a split smooth del Pezzo surface of degree $d\leq 7$ defined over $k$.
Here split means that we have $\rho(S) = \rho(S_{\overline{k}})$.
\cite{DLTT25} and \cite{Tanimoto25} studied Manin's conjecture for the trivial family $\pi: S \times \mathbb P^1 \to \mathbb P^1$. This is based on the homological sieve method developed in \cite{DLTT25}. Let us state the main result of \cite{DLTT25}:
we assume that $d = 4$. Let $\ell$ be  non-negative, rational, homogeneous, continuous, and piecewise linear function on the nef cone $\Nef_1(S)$. Let $\epsilon > 0$ be a sufficiently small rational number. Then we define the shrunken nef cone by
\[
\Nef_1(S)_\epsilon := \{ \alpha \in \Nef_1(S) \, |\, \ell(\alpha) \geq -\epsilon K_S.\alpha \}.
\]
We define the counting function by
\[
N_{\mathrm{stan}, \epsilon}(\pi, d):= \sum_{\alpha \in \Nef_1(S)_{\epsilon, \mathbb Z}, -K_S.\alpha \leq d} \#\mathrm{Mor}(\mathbb P^1, S, \alpha)(k),
\]
where $\mathrm{Mor}(\mathbb P^1, S, \alpha)$ is the morphism scheme parametrizing morphisms $s : \mathbb P^1 \to S$ such that $s_*[\mathbb P^1] = \alpha$.
Here is the main theorem of \cite{DLTT25}:

\begin{theo}[{\cite[Theorem 1.1]{DLTT25}}]
    \label{theo:dltt25}
    Fix a sufficiently small rational number $\epsilon > 0$.
    Let $q$ be a power of a prime number such that $q^\epsilon > 2^{32}$. Let $S$ be a split smooth quartic del Pezzo surface defined over $\mathbb F_q$. Then there exists a non-negative, rational, homogeneous, continuous, and piecewise linear function $\ell$ on $\Nef_1(S)$ which does not depend on $\epsilon$ and $q$ and takes positive values on a dense open cone $U \subset \Nef_1(S)$ such that
    \[
    N_{\mathrm{stan}, \epsilon}(\pi, d)\sim (1-q^{-1})^{-1}\alpha(S, \Nef_1(S)_\epsilon)\tau_{S\times \mathbb P^1}(-\mathcal K_{\mathcal X_\eta})q^d d^5,
    \]
    as $d \to \infty$.
\end{theo}

The method of the proof, developed by Das, Tosteson, and the authors, is called the homological sieve method.  Its ingredients are the following:
\begin{itemize}
    \item algebraic geometry (birational geometry of moduli spaces of rational curves on smooth quartic del Pezzo surfaces);
    \item arithmetic geometry (simplicial schemes, their homotopy theory and Grothendieck--Lefschetz trace formula);
    \item algebraic topology (the inclusion--exclusion principle and the Vassiliev type method of bar complexes); and
    \item elementary analytic number theory.
\end{itemize}

To our knowledge, this is the first time that such a geometric and topological method was used to establish Manin's conjecture over global function fields for highly non-trivial examples.
Moreover, the second author applied this method to study split del Pezzo surfaces of degree $\geq 5$ in \cite{Tanimoto25}.
Readers interested in this method should consult \cite{DLTT25}.

\begin{rema}
    Over $\mathbb C$, Cohen--Jones--Segal's conjecture (\cite[Conjecture after Definition 2]{CJS00}) roughly predicts that the homology of the irreducible components $M_{\alpha} \subset \Mor(B,X)$ stabilizes to the homology of the space of continuous maps $\mathrm{Top}(B,X)_\alpha$ as $\alpha$ increases.  \cite{DLTT25} proved a version of Cohen--Jones--Segal's conjecture for quartic del Pezzo surfaces over $\mathbb{C}$ using arguments similar to those in the proof of Theorem \ref{theo:dltt25}. This is based on the method of the bar complexes developed by Das--Tosteson in \cite{DT24} which settles Cohen--Jones--Segal's conjecture for quintic del Pezzo surfaces. An upcoming paper of Das, Tosteson, and the authors gives a general description of the relationship between Manin's conjecture and Cohen--Jones--Segal's conjecture.
\end{rema}

\bibliographystyle{alpha}
\bibliography{SRISurvey}

\end{document}